\documentclass[11pt]{article} 
\usepackage{amsmath}
\usepackage{amsfonts}
\usepackage{graphicx}
\usepackage{color}
\usepackage[latin1]{inputenc}
\usepackage[T1]{fontenc}
\usepackage[french]{babel}
\makeatletter
\def\@sect#1#2#3#4#5#6[#7]#8{%
  \ifnum #2>\c@secnumdepth
    \let\@svsec\@empty
  \else
    \refstepcounter{#1}%
    \protected@edef\@svsec{\@seccntformat{#1}\relax}%
  \fi
  \@tempskipa #5\relax
  \ifdim \@tempskipa>\z@
    \begingroup
      #6{%
        \@hangfrom{\hskip #3\relax\@svsec}%
          \interlinepenalty \@M #8\@@par}%
    \endgroup
    \csname #1mark\endcsname{#7}%
    \addcontentsline{toc}{#1}{%
      \ifnum #2>\c@secnumdepth \else
        \protect\numberline{\csname the#1\endcsname.}%
      \fi
      #7}%
  \else
    \def\@svsechd{%
      #6{\hskip #3\relax
      \@svsec #8}%
      \csname #1mark\endcsname{#7}%
      \addcontentsline{toc}{#1}{%
        \ifnum #2>\c@secnumdepth \else
          \protect\numberline{\csname the#1\endcsname.}%
        \fi
        #7}}%
  \fi
  \@xsect{#5}}
\def\@seccntformat#1{\csname the#1\endcsname.\quad}
\makeatother
\newtheorem{theo}[equation]{Th\'eor\`eme}
\newtheorem{fait}[equation]{Fait}
\newtheorem{lem}[equation]{Lemme}
\makeatletter
\renewcommand\theequation{\thesection.\arabic{equation}}
\@addtoreset{equation}{section}
\makeatother
\newcommand{\carrenoir}{\rule{0.5em}{0.5em}}
\makeatletter
\newenvironment{demo}[1][\@empty]{\textbf{D\'emonstration~%
\ifx\@empty#1:\else #1~:\fi~}}
{\hfill\carrenoir\nolinebreak\vspace{2mm}}
\makeatother
\newcommand{\oper}[2]{\newcommand{#1}{\mathop{\mathrm{#2}}\nolimits} }
\oper{\Vol}{Vol}
\newcommand{\de}{\mathrm{ d }}
\newcommand{\R}{\mathbb R}
\newcommand{\N}{\mathbb N}
\oper{\dimension}{dim}
\oper{\SO}{SO}

{\endtrivlist}%
\newenvironment{remarque}{
\refstepcounter{equation}\trivlist%
\item[\hskip \labelsep{\bfseries Remarque \theequation.\ }]}%

\DeclareSymbolFont{greek}{OML}{ptmcm}{m}{it}
\DeclareMathSymbol{\codiff}{\mathord}{greek}{"0E}
\DeclareMathSymbol{\prodint}{\mathord}{greek}{"13}

\title{Prescription du spectre du laplacien de Hodge-de~Rham dans
une classe conforme}
\author{Pierre Jammes}
\date{}
\begin{document}
\maketitle
{\small 
\textsc{Résumé.---}
Sur toute variété compacte de dimension $n\geq5$, on prescrit le volume et 
toute partie finie du spectre du laplacien de Hodge-de~Rham en restriction 
aux formes de degré $p\in[2,n-2]$, en excluant $p=n/2$ si $n$ est pair, 
et en imposant à la métrique d'appartenir à une classe conforme donnée. 
On sait que pour $n\leq4$, ainsi que pour $p=0,1,n-1,n$, et $p=n/2$ si
$n$ est pair, on ne peut pas prescrire simultanément le spectre, le volume
et la classe conforme.

Mots-clefs : formes différentielles, laplacien de Hodge-de~Rham, 
prescription de spectre, géométrie conforme.

\medskip
\textsc{Abstract.---}
 For any compact manifold of dimension $n\geq5$, we prescribe the volume
and any finite part of the spectrum of the Hodge Laplacian acting
on differential forms of degree $p\in[2,n-2]$ (exept for $p=n/2$ if $n$ is 
even), within a given conformal class. When $n\leq4$ and when $p=0,1,n-1,n$,
and $p=n/2$ if $n$ is even, this simultaneous prescription of the volume,
the spectrum and the conformal class is known to be impossible.

Keywords : differential forms, Hodge Laplacian, prescription of spectrum,
conformal geometry.

MSC2000 : 58J50, 58C40, 53A30}

\section{Introduction}

 Étant donnée une suite finie croissante de réels strictement positifs
$0<\lambda_1\leq\lambda_2\leq\ldots\leq\lambda_k$, Y.~Colin de Verdière
a montré dans \cite{cdv87} qu'on peut trouver sur toute variété
compacte de dimension supérieure ou égale à~3 une métrique riemannienne 
telle que le spectre du laplacien
agissant sur les fonctions commence par la suite $(\lambda_i)_{i=1}^k$,
et J.~Lohkamp a amélioré ce résultat en montrant dans \cite{lo96} qu'on 
pouvait prescrire le volume et certains invariants de courbure en même
temps que le spectre.  Dans le cas où la suite $(\lambda_i)_{i=1}^k$ 
est strictement croissante, c'est-à-dire en supposant que les valeurs 
propres prescrites sont simples, des résultats du même type ont été obtenus
pour d'autres opérateurs: M.~Dahl prescrit dans \cite{da05} le début 
du spectre de l'opérateur de Dirac sur les variétés compactes et
P.~Guérini a montré dans \cite{gu04} que sur les variétés compactes et 
les domaines euclidiens, on peut prescrire simultanément le volume et 
toute partie finie du spectre du laplacien de Hodge-de~Rham, qui agit 
sur les formes différentielles.

Notre but est de montrer, dans le cas du laplacien de Hodge-de~Rham,
que si on munit une variété compacte $M^n$ de dimension $n$ d'une métrique
riemannienne $g$ quelconque, on peut obtenir le volume et la partie finie 
du spectre (sans multiplicité)
souhaités en effectuant uniquement des déformations conformes à
partir de $g$. Il faut noter qu'une telle prescription simultanée
du spectre, du volume et de la classe conforme est spécifique aux formes
différentielles. Elle est impossible pour le laplacien agissant sur
les fonctions: si on fixe le volume et la classe conforme sur une variété
compacte donnée, on ne peut pas
rendre les valeurs propres arbitrairement grandes (\cite{esi86}, 
\cite{ko93}). Sous les mêmes contraintes, on ne peut pas rendre
les valeurs propres non nulles de l'opérateur de Dirac arbitrairement
petites (\cite{lo86}, \cite{am03}). Dans le cas des formes différentielles,
on se heurte au problème que le spectre des $1$-formes contient le 
spectre des fonctions, et pour les formes de degré $n/2$ quand $n$ est
pair il y a une obstruction du même type que pour l'opérateur de Dirac
(cf. \cite{ja06}, ainsi que la remarque \ref{intro:rqsob} ci-dessous).
On va montrer qu'on peut prescrire le spectre pour les autres degrés. 
La comparaison
de ces différents résultats permet de mesurer la rigidité qu'apporte
le fait de fixer le volume et la classe conforme: on frôle les limites
des possibilités de prescription.

Précisons quelques notations: Si $(M^n,g)$ est une variété riemannienne
compacte orientable de dimension $n$, le laplacien $\Delta^p$ agissant sur
l'espace $\Omega^p(M)$ des $p$-formes différentielles est défini par
$\Delta=\de\codiff+\codiff\de$ où $\codiff$ désigne la codifférentielle, 
et son spectre sera noté 
\begin{equation}
0=\lambda_{p,0}(M,g)<\lambda_{p,1}(M,g)\leq\lambda_{p,2}(M,g)\leq\ldots
\end{equation}
où les valeurs propres non nulles sont répétées s'il y a multiplicité.
La multiplicité de la valeur propre nulle, si elle existe, est un invariant
topologique: c'est le nombre de Betti $b_p(M)$. 

L'espace des $p$-formes coexactes est stable par le laplacien, et
on notera 
\begin{equation}
0<\mu_{p,1}(M,g)\leq\mu_{p,2}(M,g)\leq\ldots
\end{equation} 
le spectre du laplacien restreint à cet espace. 
Par théorie de Hodge, le spectre 
$(\lambda_{p,i}(M,g))_{i\geq1}$ est la réunion de $(\mu_{p,i}(M,g))_i$ et 
$(\mu_{p-1,i}(M,g))_i$. et on a de plus
$\mu_{p,i}(M,g)=\mu_{n-p-1,i}(M,g)$ pour tout $p$ et $i$ si $M$ n'a pas
de bord.
Le spectre complet du laplacien se déduit alors des $\mu_{p,i}(M,g)$ pour
$p\leq\frac{n-1}2$, ce sont donc ces valeurs propres qu'on va chercher à 
prescrire. 
On exclut le cas $p=0$ puisque $(\mu_{0,i}(M,g))$
est le spectre des fonctions, pour lequel la prescription est
impossible comme on l'a déjà remarqué. On écarte aussi le cas 
$p=\left[\frac{n-1}2\right]$ pour lequel on ne peut pas rendre les
valeurs propres arbitrairement petites (voir remarque \ref{intro:rqsob}).
En particulier, si la dimension de $M$ vérifie $n\leq4$, la prescription 
dans une classe
conforme et à volume fixé d'une valeur propre quelconque est impossible,
quel que soit le degré.
Compte tenu de ces remarques, en supposant que $n\geq5$ et en notant $k$ 
l'entier tel que $n=2k+3$ ou 
$2k+4$, on va prescrire toute partie finie des $\mu_{p,i}(M,g)$ pour 
$1\leq p\leq k$.

\begin{theo}\label{intro:th}
Soit $M$ une variété compacte, connexe, orientable et sans bord de dimension 
$n=2k+3$ ou $2k+4$ où $k\in\N^*$,
$C$ une classe conforme de métriques riemanniennes sur $M$,
$V_0$ un réel strictement positif et $N\geq1$ un entier. 
On se donne pour 
tout entier $p\in\{1,\ldots,k\}$ une suite de réels 
$0<\nu_{p,1}<\nu_{p,2}<\ldots<\nu_{p,N}$.

Il existe une métrique $g\in C$ telle que
\begin{itemize}
\item $\mu_{p,i}(M,g)=\nu_{p,i}$ pour tout $i\leq N$ et $p\in\{1,\ldots,k\}$;
\item $\mu_{k+1,1}(M,g)>\sup_{p,i}\{\nu_{p,i}\}$;
\item $\Vol(M,g)=V_0$.
\end{itemize}
\end{theo}
\begin{remarque}
La minoration $\mu_{k+1,1}(M,g)>\sup_{p,i}\{\nu_{p,i}\}$ assure qu'on a
l'égalité
$\lambda_{k+1,i}(M,g)=\mu_{k,i}(M,g)$ pour $i\leq N$. On peut donc
prescrire les $N$ premières valeurs propres des $(k+1)$-formes, les formes 
propres correspondantes étant alors exactes, de valeurs propres égales
à $(\mu_{k,i}(M,g))_{i=1}^N$. Si $n$ est impair,
on prescrit ainsi le spectre en tout degré $2\leq p\leq n-2$. En
dimension paire, le degré $p=n/2=k+2$ fait exception. En degré $1$ et $n-1$
on ne prescrit pas arbitrairement le début du spectre car on ne contrôle 
pas les $\mu_{0,i}(M,g)$, mais on peut assurer que les valeurs 
$\nu_{1,1},\ldots,\nu_{1,N}$ sont contenues dans 
$(\lambda_{1,i}(M,g))_{i\geq1}$ et $(\lambda_{n-1,i}(M,g))_{i\geq1}$.
\end{remarque}
\begin{remarque} 
Les valeurs propres $\lambda_{p,i}(M,g)$ prescrites
sont simples, ou de multiplicité 2 si on fait en sorte que $\nu_{p,i}=
\nu_{p-1,j}$ pour des valeurs quelconques de $i$ et $j$. Le problème
de prescrire arbitrairement la multiplicité ne serait-ce que d'une
valeur propre reste à notre connaissance ouvert, tant pour le laplacien
de Hodge-de~Rham que pour l'opérateur de Dirac.
\end{remarque}

Une étape clef de la démonstration du théorème \ref{intro:th} consiste
à montrer que sur la sphère, on peut prescrire une valeur propre,
toutes les autres valeurs propres étant arbitrairement grandes,
le volume étant majoré et la classe conforme étant fixée. On va montrer un 
résultat équivalent, à savoir qu'on peut faire tendre une valeur propre 
non nulle
vers zéro en déformant la sphère de manière conforme, les autres valeurs
propres étant minorées et le volume étant fixé:
\begin{lem}\label{intro:lem}
Soit $n\geq5$ un entier, $k$ l'entier tel que $n=2k+3$ ou $n=2k+4$ et $C$ 
une classe conforme sur $S^n$. Pour tout réel $V>0$ et tout entier $1\leq
p\leq k$
il existe une famille de métriques $(g_\varepsilon)_{0<\varepsilon<1}$ 
contenue dans $C$ et une constante $c>0$ telles que 
$\mu_{p,1}(S^n,g_\varepsilon)<\varepsilon$, $\mu_{p,2}(S^n,g_\varepsilon)>c$, 
$\mu_{q,1}(S^n,g_\varepsilon)>c$ pour $1\leq q\leq k+1$, $q\neq p$, et 
$\Vol(S^n,g_\varepsilon)=V$.
\end{lem}
\begin{remarque}\label{intro:rqsob}
Ce lemme ne se généralise pas aux formes différentielles de degré $k+1$; 
on montre en effet dans \cite{ja06} qu'une inégalité de Sobolev permet 
de minorer
uniformément $\mu_{k+1,1}(M,g)\Vol(M,g)^{2/n}$ sur une classe conforme 
par une constante strictement positive. C'est la raison pour laquelle 
on ne peut pas prescrire $\mu_{k+1,i}(M,g)$ dans le théorème \ref{intro:th}.
\end{remarque}
\begin{remarque} Le problème de faire tendre des valeurs propres vers 0
dans une classe conforme à volume fixé avait été posé par B.~Colbois dans
\cite{co04} et était resté ouvert. Le lemme~\ref{intro:lem} y répond 
partiellement et la technique utilisée permet d'obtenir un grand nombre
de petite valeurs propres sur une variété quelconque (voir 
remarque~\ref{pvp:rq2}). Ce problème est totalement résolu par le
théorème~\ref{intro:th} et les résultats de \cite{ja06}.
\end{remarque}

\section{Petites valeurs propres dans une classe conforme}\label{pvp}
\subsection{Quasi-isométries et extrema conformes du spectre}\label{pvp:dodziuk}
 Commençons par rappeler le lemme suivant, dû à J.~Dodziuk, qui permet
de comparer les spectres de deux métriques dont on connaît le rapport
de quasi-isométrie et que nous utiliserons à plusieurs reprises au
cours des démonstrations du théorème \ref{intro:th} et du lemme 
\ref{intro:lem} :
\begin{lem}[\cite{do82}]\label{pvp:lem}
Soit $g$ et $\tilde g$ deux métriques riemanniennes sur une variété 
compacte $M$ de dimension $n$, et $\tau$ une constante strictement positive. 
Si les deux métriques vérifient $\frac1\tau g\leq\tilde g\leq\tau g$, alors
$$\frac1{\tau^{3n-1}}\lambda_{p,k}(M,g)\leq\lambda_{p,k}(M,\tilde g)
\leq\tau^{3n-1}\lambda_{p,k}(M,g),$$
pour tout entiers $k\geq0$ et $p\in[0,n]$.
\end{lem}

Une première conséquence du lemme de Dodziuk est qu'il suffit de démontrer
le lemme \ref{intro:lem} pour une classe conforme particulière, le 
résultat général s'en déduira : supposons que le lemme soit vrai pour
une classe conforme $C$, et donnons-nous une autre classe conforme $C'$,
ainsi que deux métriques $g\in C$ et $g'\in C'$, un réel $V>0$ 
et un entier $p\in[1,k]$. Le lemme \ref{intro:lem} nous dit 
que pour tout $\varepsilon>0$ et tout $\tau>1$, il existe une fonction
$h\in C^\infty(M)$ strictement positive telle que 
$\mu_{p,1}(M,h^2g)<\tau^{-1-3n}\varepsilon$, $\mu_{p,2}(M,h^2g)>\tau^{3n+1}c$,
$\mu_{q,1}(M,h^2g)>\tau^{3n+1}c$ pour tout $0<p\leq k+1$, $q\neq p$, 
et $\Vol(M,h^2g)=V$.

Par compacité de $M$, il existe une constante $\tau>1$ telle que 
$\frac1\tau g\leq g'\leq\tau g$. 
Or, les métriques $h^2g$ et $h^2g'$ sont
liées par le même rapport de quasi-isométrie que $g$ et $g'$, c'est-à-dire
que $\frac1\tau h^2g\leq h^2g'\leq\tau h^2g.$ 
On en déduit alors que 
$\mu_{p,1}(M,h^2g')<\tau^{-2}\varepsilon$, $\mu_{p,2}(M,h^2g')>\tau^2c$,
$\mu_{q,1}(M,h^2g')>\tau^2c$ pour $q\neq p$
et $\tau^{-n}V\leq\Vol(M,h^2g)\leq\tau^nV$. Après renormalisation du
volume par homothétie, on a $\mu_{p,1}(M,h^2g')<\varepsilon$,
$\mu_{p,2}(M,h^2g')>c$, $\mu_{q,1}(M,h^2g')>c$ et $\Vol(M,h^2g)=V$, 
et donc le lemme \ref{intro:lem} est vrai pour la classe conforme $C'$.

\subsection{Construction d'une petite valeur propre sur la sphère}
On va maintenant aborder la démonstration du lemme~\ref{intro:lem}, qui
se déroule en deux étapes. Dans un premier temps on va montrer qu'on
peut obtenir une petite valeur propre de degré $p\in[1,k]$ dans une
classe conforme donnée et à volume fixé. On vérifiera ensuite que c'est
bien la seule petite valeur propre.

Le principe de la construction de cette petite valeur propre est le suivant:
on plonge une sphère $S^p$ dans $S^n$ et on écrase la métrique en
dehors d'un voisinage tubulaire de $S^p$. On peut alors choisir une
forme test qui a un petit quotient de Rayleigh en prolongeant la forme
volume de $S^p$.

On fixe donc un entier $p\in[1,k]$, et on considère un plongement
$i:S^p\times B^{n-p}\hookrightarrow S^n$ où $B^{n-p}$ est la boule
de dimension $n-p$, l'image de $i$ étant un voisinage tubulaire d'une 
sous-variété de $S^n$ difféomorphe à $S^p$. Dans la suite, on identifiera
$S^p\times B^{n-p}$ avec son image par $i$ qu'on notera $\Omega$.
On identifiera aussi les formes volumes $\de v_{S^p}$ et $\de v_{B^{n-p}}$
de $S^p$ et $B^{n-p}$ avec leur relevé sur $\Omega$.

Comme on l'a montré au paragraphe~\ref{pvp:dodziuk}, il suffit de montrer 
le lemme pour
la classe conforme d'un métrique particulière de $S^n$ qu'on choisit comme 
suit : on muni $B^{n-p}$ d'une
métrique euclidienne de rayon $R$, $S^p$ d'une métrique quelconque, 
$\Omega$ de la métrique produit associée, et $S^n$ d'une métrique $g$
qui prolonge la métrique sur $\Omega$. On choisit $R$ suffisamment
petit pour qu'on puisse choisir $g$ telle que $\Vol(S^n,g)=V$, et on va 
montrer qu'on peut trouver une métrique $g_\varepsilon$ conforme à $g$ telle 
que $\mu_{p,1}(S^n,g_\varepsilon)<\varepsilon$ et $\Vol(S^n,g_\varepsilon)=V$.
 
On note $r$ la coordonnée radiale sur $B^{n-p}$ et on se donne une 
fonction continue $f$ sur $B^{n-p}$ qui ne dépend que de $r$ telle que
$f(0)=1$, $f(R)=0$. On construit une $p$-forme test $\omega$ sur $\Omega$
en posant $\omega=f\de v_{S^p}$, et on l'étend en une $p$-forme sur $S^n$ 
par $\omega=0$ en dehors de $\Omega$. Sur $\Omega$, la forme $\omega$ vérifie 
$\de\omega=f'\de r\wedge\de v_S$ et $(-1)^{n(p+1)+1}\codiff\omega=
*\de(fh^{n-2p}\de v_B)=0$ car $\de r\wedge \de v_B=0$.

\begin{figure}[h]
\begin{center}
\begin{picture}(0,0)%
\includegraphics{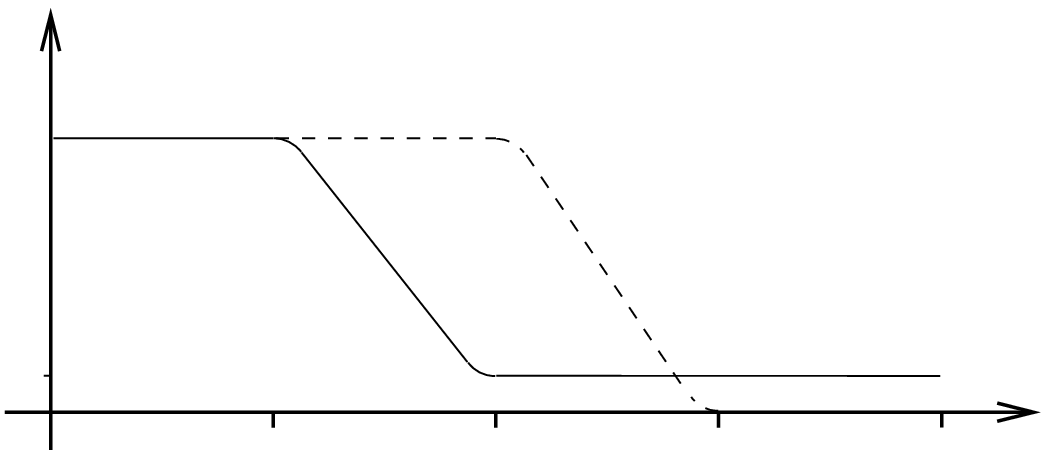}%
\end{picture}%
\setlength{\unitlength}{4144sp}%
\begingroup\makeatletter\ifx\SetFigFontNFSS\undefined%
\gdef\SetFigFontNFSS#1#2#3#4#5{%
  \reset@font\fontsize{#1}{#2pt}%
  \fontfamily{#3}\fontseries{#4}\fontshape{#5}%
  \selectfont}%
\fi\endgroup%
\begin{picture}(4781,2163)(3522,-2666)
\put(3610,-1162){\makebox(0,0)[lb]{\smash{{\SetFigFontNFSS{8}{9.6}{\rmdefault}{\mddefault}{\updefault}{\color[rgb]{0,0,0}$1$}%
}}}}
\put(3610,-2498){\makebox(0,0)[lb]{\smash{{\SetFigFontNFSS{8}{9.6}{\rmdefault}{\mddefault}{\updefault}{\color[rgb]{0,0,0}$0$}%
}}}}
\put(5569,-1349){\makebox(0,0)[lb]{\smash{{\SetFigFontNFSS{8}{9.6}{\rmdefault}{\mddefault}{\updefault}{\color[rgb]{0,0,0}$f$}%
}}}}
\put(5009,-1767){\makebox(0,0)[lb]{\smash{{\SetFigFontNFSS{8}{9.6}{\rmdefault}{\mddefault}{\updefault}{\color[rgb]{0,0,0}$h$}%
}}}}
\put(4692,-2615){\makebox(0,0)[lb]{\smash{{\SetFigFontNFSS{8}{9.6}{\rmdefault}{\mddefault}{\updefault}{\color[rgb]{0,0,0}$\frac{R}4$}%
}}}}
\put(5721,-2615){\makebox(0,0)[lb]{\smash{{\SetFigFontNFSS{8}{9.6}{\rmdefault}{\mddefault}{\updefault}{\color[rgb]{0,0,0}$\frac{R}2$}%
}}}}
\put(6711,-2615){\makebox(0,0)[lb]{\smash{{\SetFigFontNFSS{8}{9.6}{\rmdefault}{\mddefault}{\updefault}{\color[rgb]{0,0,0}$\frac{3R}4$}%
}}}}
\put(4138,-2581){\makebox(0,0)[lb]{\smash{{\SetFigFontNFSS{8}{9.6}{\rmdefault}{\mddefault}{\updefault}{\color[rgb]{0,0,0}$I_1$}%
}}}}
\put(5156,-2581){\makebox(0,0)[lb]{\smash{{\SetFigFontNFSS{8}{9.6}{\rmdefault}{\mddefault}{\updefault}{\color[rgb]{0,0,0}$I_2$}%
}}}}
\put(6174,-2581){\makebox(0,0)[lb]{\smash{{\SetFigFontNFSS{8}{9.6}{\rmdefault}{\mddefault}{\updefault}{\color[rgb]{0,0,0}$I_3$}%
}}}}
\put(7231,-2581){\makebox(0,0)[lb]{\smash{{\SetFigFontNFSS{8}{9.6}{\rmdefault}{\mddefault}{\updefault}{\color[rgb]{0,0,0}$I_4$}%
}}}}
\put(7757,-2610){\makebox(0,0)[lb]{\smash{{\SetFigFontNFSS{8}{9.6}{\rmdefault}{\mddefault}{\updefault}{\color[rgb]{0,0,0}$R$}%
}}}}
\put(3599,-2240){\makebox(0,0)[lb]{\smash{{\SetFigFontNFSS{8}{9.6}{\rmdefault}{\mddefault}{\updefault}{\color[rgb]{0,0,0}$\eta$}%
}}}}
\end{picture}%
\end{center}
\caption{\label{pvp:graphe}}
\end{figure}

Précisons le choix de $f$ et de la nouvelle métrique. On veut contracter
la métrique dans un domaine de $S^n$ qui contient le support de $f'$. Étant 
donné un réel $0<\eta\ll1$, on pose $\bar g=h^2g$ où $h$ est une fonction qui, 
sur $S^p\times B^{n-p}$, ne dépend que de $r$ et vérifie $h(0)=1$, $h(1)=\eta$ 
et qui vaut $\eta$ en dehors de $S^p\times B^{n-p}$.
On découpe l'intervalle $[0,R]$ en quatre intervalles $I_i=[\frac{(i-1)R}4,
\frac{iR}4]$ et on pose :
\begin{itemize}
\item $h(r)=1$ sur $I_1$;
\item $h(r)=\eta$ sur $I_3$, $I_4$;
\item $f(r)=1$ sur $I_1$ et $I_2$;
\item $f(r)=0$ sur $I_4$.
\end{itemize}
On prolonge $h$ sur $I_2$ et $f$ sur $I_3$ de manière lisse et 
monotone (voir figure~\ref{pvp:graphe}).

Le quotient de Rayleigh de la forme $\omega$ pour la métrique $\bar g$ est

\begin{eqnarray}
R(\omega)&=&\frac{\|\de\omega\|_{\bar g}^2}{\|\omega\|_{\bar g}^2}=
\frac{\displaystyle\int_{S^n} h^{n-2p-2}|\de\omega|_g^2\de v_g}
{\displaystyle\int_{S^n} h^{n-2p}|\omega|_g^2\de v_g}\\\nonumber
&=& \frac{\displaystyle\int_{S^n} h^{n-2p-2}|f'\de r\wedge\de v_{S^p}|^2\
\de v_g} {\displaystyle\int_{S^n} h^{n-2p}|f|^2\de v_g}.
\end{eqnarray}
 La fonction $f'$ est nulle en dehors de $I_3$, par conséquent 
\begin{equation}\label{pvp:eq1}
R(\omega)=\frac{\displaystyle\eta^{n-2p-2}\left(\int_{r\in I_3}|f'|^2\de v_g
\right)}{\displaystyle\int_{S^n} h^{n-2p}|f|^2\de v_g}\leq
\eta^{n-2p-2}\frac{\displaystyle\int_{r\in I_3}|f'|^2\de v_g}
{\displaystyle \int_{r\in I_1}\de v_g}.
\end{equation}
 Dans l'inégalité (\ref{pvp:eq1}), le quotient du membre de gauche est
uniformément majoré par rapport à $\eta$.  On a supposé que $p\leq k$, 
donc $n-2p-2$ est strictement positif. En choisissant $\eta$ suffisamment 
petit, on obtient une fonction $h_\varepsilon$ et une métrique 
$\bar g_\varepsilon=h_\varepsilon^2g$ telle que $R(\omega)<\varepsilon$, 
et par conséquent $\mu_{p,1}(M,\bar g_\varepsilon)<\varepsilon$ car la 
forme test est cofermée et la sphère $S^n$ n'a pas de cohomologie 
en degré $p$. Le volume de $S^n$ ne fait que diminuer lorsque $\eta$ tend 
vers zéro. Si on le normalise par une homothétie on obtient une métrique 
$g_\varepsilon$ telle que $\Vol(S^n,g_\varepsilon)=V$, la petite valeur 
propre restant inférieure à $\varepsilon$, ce qui conclut la première partie
de la démonstration du lemme~\ref{intro:lem}.

\begin{remarque}\label{pvp:rq1}
 En dehors du domaines $\Omega$, la métrique ne subit qu'une 
homothétie. Cette propriété simplifiera la prescription du spectre.
\end{remarque}
\begin{remarque}\label{pvp:rq2}
 Le procédé que nous venons de décrire se généralise aisément pour
construire un grand nombre de petites valeurs propres pour plusieurs
degrés simultanément sur une variété compacte $M$ quelconque : à partir 
d'un nombre arbitraire de sphères plongées dans $M$, on peut construire 
autant de formes test dont le quotient de Rayleigh tend vers zéro quand
on écrase la métrique en dehors de voisinages tubulaires tous disjoints de
ces sphères. Cependant, si la variété a des nombres de Betti non nuls,
on ne contrôle pas précisément le nombre de petites valeurs propres non 
nulles obtenu.
\end{remarque}

\subsection{Contrôle du nombre de petites valeurs propres}

Pour achever la démonstration du lemme \ref{intro:lem}, on doit s'assurer
que la petite valeur propre construite au paragraphe précédent est
bien la seule. Nous allons pour cela faire appel à un lemme dû à J.~McGowan
qui permet, étant donné un recouvrement de la variété par des
ouverts à bords lisses, de minorer une partie du spectre de la variété
en fonction du spectre des ouverts du recouvrement et de leurs intersections.
Sur les domaines de la variété, le spectre que nous considérerons sera 
toujours celui du laplacien agissant sur les formes vérifiant la condition 
de bord absolue, qui généralise la condition de Neumann, à savoir :
\begin{equation}
\left\{\begin{array}{l}
j^*\prodint_\nu\omega=0\\
j^*\prodint_\nu\de\omega=0,\end{array}\right.
\end{equation}
où $j:\partial M\hookrightarrow M$ désigne l'inclusion canonique, $\nu$ la
normale au bord et $\prodint_\nu$ le produit intérieur par $\nu$.

\begin{lem}\label{pvp:mg}
Soit $(M,g)$ une variété compacte de dimension $n$, $(U_i)_{i=1}^K$ un 
recouvrement de $M$ par des ouverts n'ayant pas d'intersections d'ordre
supérieur ou égal à~3 et $\rho_i$ une partition
de l'unité relative à $(U_i)_{i=1}^K$. Il existe des constantes $a,b>0$
ne dépendant que de $n$ telles que
$$\mu_{q,k_q}(M,g)\geq\frac a{\displaystyle\sum_{i=1}^K
\left(\frac1{\mu_i}+\displaystyle\sum_{U_i\cap U_j\neq\emptyset}
\left(\frac{b\cdot c_\rho}{\mu_{ij}}+1\right)
\left(\frac1{\mu_i}+\frac1{\mu_j}\right)\right)}$$
avec $k_q=1+\sum_{i,j}\dimension\mathcal H^q(\bar U_i\cap\bar U_j)$ où 
$\mathcal H^q$ désigne l'espace des $p$-formes harmoniques avec condition 
de bord absolue, $\mu_i=\mu_{q,1}(\bar U_i)$, $\mu_{ij}=\mu_{q-1,1}
(\bar U_i\cap\bar U_j)$ et $c_\rho=\sup_i\|\nabla\rho_i\|_\infty^2$.
\end{lem}
Ce lemme a été démontré dans \cite{mc93} pour les $1$-formes pour
un recouvrement fini quelconque, et G.~Gentile et V.~Pagliara ont remarqué
dans \cite{gp95} qu'il se généralise aux formes de degré quelconque si
on suppose que les ouverts du recouvrement n'ont pas d'intersection d'ordre
supérieur ou égal à~3.

Nous allons appliquer ce lemme au recouvrement constitué des deux ouverts
$U_1$ formé de l'intérieur de $\Omega$ et $U_2$ formé de $S^n$ privé 
des points de $\Omega$ tels que $r\leq \frac{3R}4$.  Par commodité,
on appliquera le lemme \ref{pvp:mg} pour la métrique $\bar g_\varepsilon$
au lieu de $g_\varepsilon$. Comme $\Vol(S^n,\bar g_\varepsilon)$ reste
uniformément minoré quand $\varepsilon\to0$, cela ne change pas
significativement le résultat.

 Le nombre $k_q$ du lemme~\ref{pvp:mg} vaut $1$ pour $q=p$
et 0 pour les autres degrés $q\leq k+1$. Par conséquent, la minoration
des termes $\mu_1$, $\mu_2$ et $\frac{\mu_{12}}{c_\rho}$ du lemme~\ref{pvp:mg}
permettrait de conclure la démonstration. 

On peut noter
qu'au cours de la déformation conforme créant la petite valeur propre,
l'intersection $U_1\cap U_2$ ne subit qu'une homothétie.
En choisissant une partition de l'unité indépendante de $\bar g_\varepsilon$, 
le rapport $\frac{c_\rho}{\mu_{12}}$ dans le lemme \ref{pvp:mg}
est invariant par homothétie, il est donc indépendant de $\varepsilon$.
Par ailleurs, la métrique sur $U_2$ subit une homothétie qui fait tendre 
son spectre non nul vers l'infini. 

Pour minorer $\mu_{q,k_q}(S^n,\bar g_\varepsilon)$ à l'aide du 
lemme~\ref{pvp:mg} il reste donc à
minorer le spectre de $U_1$, pour lequel un difficulté apparaît :
la métrique sur ce domaine n'est pas fixe puisqu'on l'écrase sur
un voisinage du bord, et il n'est pas clair \emph{a priori} que cette
déformation ne produit pas de petites valeurs propres : on a vu précédemment
que si une telle déformation se produit ailleurs que près du bord, elle peut
effectivement faire tendre une valeur propres vers zéro.

\begin{fait}
Il existe une constante $c>0$ indépendante du choix de $h_\varepsilon$ 
telle que $\mu_{q,1}(\Omega,h_\varepsilon^2g)>c$, pour tout $1\leq q\leq k+1$.
\end{fait}
\begin{demo} 
On utilise une nouvelle fois le lemme de Dodziuk pour se ramener à une
classe conforme particulière en identifiant $B^{n-p}$ à la réunion d'un 
cylindre $[0,1]\times S^{n-p-1}$ (muni d'une métrique produit) et d'un
hémisphère (voir figure~\ref{pvp:mgfig}), en munissant $S^p$ de sa métrique 
canonique et $\Omega$ de la métrique produit.

\begin{figure}[h]
\begin{center}
\begin{picture}(0,0)%
\includegraphics{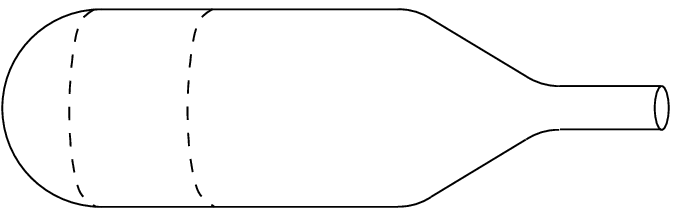}%
\end{picture}%
\setlength{\unitlength}{4144sp}%
\begingroup\makeatletter\ifx\SetFigFontNFSS\undefined%
\gdef\SetFigFontNFSS#1#2#3#4#5{%
  \reset@font\fontsize{#1}{#2pt}%
  \fontfamily{#3}\fontseries{#4}\fontshape{#5}%
  \selectfont}%
\fi\endgroup%
\begin{picture}(3065,1729)(1116,-1354)
\put(1176,-1001){\makebox(0,0)[lb]{\smash{{\SetFigFontNFSS{12}{14.4}{\rmdefault}{\mddefault}{\updefault}{\color[rgb]{0,0,0}$\underbrace{\hspace{2cm}}$}%
}}}}
\put(1531, 19){\makebox(0,0)[lb]{\smash{{\SetFigFontNFSS{12}{14.4}{\rmdefault}{\mddefault}{\updefault}{\color[rgb]{0,0,0}$\overbrace{\hspace{5.7cm}}$}%
}}}}
\put(2721,204){\makebox(0,0)[lb]{\smash{{\SetFigFontNFSS{12}{14.4}{\rmdefault}{\mddefault}{\updefault}{\color[rgb]{0,0,0}$\Omega_2$}%
}}}}
\put(1556,-1281){\makebox(0,0)[lb]{\smash{{\SetFigFontNFSS{12}{14.4}{\rmdefault}{\mddefault}{\updefault}{\color[rgb]{0,0,0}$\Omega_1$}%
}}}}
\end{picture}%
\end{center}
\caption{\label{pvp:mgfig}}
\end{figure}

On recouvre $\Omega$ par deux domaines, $\Omega_1$ défini comme
le produit de $S^p$ et de réunion de l'hémisphère et de $[0,\frac14]\times
S^{n-p-1}$, et $\Omega_2$ formé de $[0,1]\times S^{n-p-1}\times S^p$.
On va minorer le spectre sur chacun des deux domaines et en déduire une
minoration sur $\Omega$. On ne peut pas appliquer le lemme de 
McGowan pour minorer la totalité du spectre non nul car la cohomologie des 
intersections des domaines est non triviale, mais dans cette situation on peut
en adapter la démonstration pour obtenir un résultat plus précis.

Remarquons d'abord que la métrique sur $\Omega_1$ est fixe, son spectre
est donc uniformément minoré par rapport à $\varepsilon$. On doit
ensuite minorer le spectre de $\Omega_2$. Sur ce domaine, la métrique
$g$ est le produit des métriques canoniques de $[0,1]$, $S^p$ et $S^{n-p-1}$.
Elle est donc invariante sous l'action des groupes d'isométries
$\SO(p+1)$ et $\SO(n-p)$ des sphères $S^p$ et $S^{n-p-1}$, et il en est
de même pour la métrique $h_\varepsilon^2g$.
Or, on peut montrer (voir \cite{ja04}, théorème~1.18 et section~4) que
lorsque le cercle agit par isométrie sur une variété compacte, il existe
une constante dépendant uniquement de la longueur maximale des orbites telle
que si une valeur propre est inférieure à cette constante, les formes
propres correspondantes sont invariantes sous l'action du cercle. Le groupe 
$\SO(p+1)\times\SO(n-p)$ étant engendré par des cercles dont la longueur
des orbites est uniformément majorée par rapport à $\varepsilon$,
on peut se restreindre aux formes invariantes par cette action. Remarquons
en outre que les seules formes invariantes de la sphère sont les fonctions
constantes et les multiples de la forme volume canonique.

Ces remarques nous permettent de nous ramener à un problème unidimensionnel:
une forme propre de $(\Omega_2,h^2g)$ de degré compris entre 1 et $k+1$ et 
dont la valeur propre est petite est nécessairement de la forme 
\begin{equation}
\omega=f(t)\de v_{S^p} \textrm{ ou } \omega=f(t)\de v_{S^p}\wedge\de t,
\end{equation}
où $f$ est une fonction sur $[0,1]$ et $\de v_{S^p}$ la forme volume
de $S^p$, qu'on identifie avec leur relevé à $\Omega_2$. L'image 
de $\omega$ par $\de$ étant aussi invariante, on peut
affirmer que  $f(t)\de v_{S^p}\wedge\de t$ est fermée. Comme il suffit de
minorer le spectre des formes coexactes, on va raisonner sur
$\omega=f(t)\de v_{S^p}$.

On va calculer explicitement le laplacien de $\omega$ en fonction de $f$ et
$h$. Remarquons d'abord que $\de*\omega=\de(h(t)^{n-2p}f(t)
\de v_{S^{n-p-1}}\wedge\de t)=0$, c'est-à-dire que $\codiff\omega=0$. 
On a par ailleurs $\de\omega=f'(t)\de t\wedge\de v_{S^p}$, donc
\begin{eqnarray}
\de*\de\omega&=&(-1)^{n-1}\de(f'h^{n-2p-2})\de v_{S^{n-p-1}}\\
&=&(-1)^l[h^{n-2p-2}f''+(n-2p-2)f'h'h^{n-2p-3}]\de v_{S^{n-p-1}}\wedge\de t
\nonumber
\end{eqnarray}
et finalement 
\begin{equation}
\Delta\omega=\codiff\de\omega=-[h^{-2}f''+(n-2p-2)f'h'h^{-3}]\de v_{S^p}.
\end{equation}
 Le même calcul montre que sur la variété $S^{n-2p-1}\times[0,1]$
munie de la métrique $h^2g'$ où $g'$ est la métrique produit, la fonction
$f$ vérifie $\Delta f=-(h^{-2}f''+(n-2p-2)f'h'h^{-3})$. Si $\omega$ est
une forme propre de valeur propre $\lambda$, alors $f$ est une fonction 
propre sur $(S^{n-2p-1}\times[0,1],h^2g')$ de même valeur propre.
Or, si on choisit la famille
$h_\varepsilon$ de sorte qu'elle tende simplement vers la fonction
caractéristique de $[0,1/2]$, on sait que cette déformation ne produit pas
de petites valeurs propres pour les fonctions (c'est un corollaire immédiat
du théorème~III.3 de \cite{cdv86}). Le spectre de $\Omega_2$ est donc
bien uniformément minoré.

On doit maintenant déduire une minoration du spectre de $\Omega$ des
minorations des spectres de $\Omega_1$ et $\Omega_2$. Pour les degrés
autres que $p$ on peut appliquer le lemme~\ref{pvp:mg}. En degré $p$ 
ce lemme ne minore que la 2\ieme{} valeur propre, on va voir comment
améliorer ce résultat. 
On se référera implicitement au début de la section 2 de l'article 
\cite{mc93} de J.~McGowan pour les notions de théorie de Hodge et de 
théorie spectrale des variétés à bord que nous utiliserons. Rappelons-en
deux points techniques: pour minorer $\mu_{p,1}$ il suffit de trouver une
constante $C>0$ telle que pour toute $(p+1)$-forme exacte $\phi$, il existe
une $p$-forme $\psi$ telle que $\de\psi=\phi$ et 
$\|\psi\|/\|\phi\|\leq C$, et à l'exception 
des formes harmoniques, on a pas besoin de se restreindre aux formes 
différentielles vérifiant les conditions de bord.

Soit $\phi$ une $(p+1)$-forme exacte sur $\Omega$. Pour tout $i=1,2$, on note
$\phi_i$ sa restriction au domaine $\Omega_i$. Chaque $\phi_i$ est
une forme exacte, donc il existe sur chacun des $\Omega_i$ 
une forme $\psi_i$ telle que $\de\psi_i=\phi_i$, et on peut choisir
$\psi_i$ telle que $\mu_{p,1}(\Omega_i)\|\psi_i\|^2\leq\|\phi_i\|^2$.

Sur l'intersection $\Omega_{12}=\Omega_1\cap\Omega_2$, on peut définir la forme
$\omega=\psi_2-\psi_1$. Elle vérifie $\de\omega=\de\psi_2-\de\psi_1=0$,
on peut donc l'écrire sous la forme $\omega=\alpha+\de\beta$, où $\alpha$ 
est une $p$-forme harmonique ---~avec conditions de bord absolues~--- de 
$\Omega_{12}\simeq S^p\times S^{n-p-1}\times [0,\frac14]$ et qui est donc 
nécessairement proportionnelle au relevé de $\de v_{S^p}$.  Cette forme 
s'étend naturellement en une forme harmonique (pour la métrique $g$) de 
$\Omega_2\simeq S^p\times S^{n-p-1}\times [0,1]$, qui est de longueur 
constante pour $g$ et qu'on notera encore $\alpha$.

Soit $(\rho_i)$ une partition de l'unité pour le recouvrement $(\Omega_i)$.
La forme $\rho_1\beta$ (resp. $\rho_2\beta$) qui est définie sur 
$\Omega_{12}$ se prolonge naturellement par 0 en dehors de $\Omega_{12}$ pour 
donner une forme sur $\Omega_2$ (resp. $\Omega_1$).
On définit alors sur les domaines $\Omega_i$ les formes 
$\bar\psi_i$ par
\begin{equation}
\bar\psi_1=\psi_1+\de(\rho_2\beta)\textrm{ et }
\bar\psi_2=\psi_2-\alpha-\de(\rho_1\beta).
\end{equation}
Ces formes vérifient $\de\bar\psi_i=\de\psi_i=\phi_i$ pour $i=1,2$, et sur 
$\Omega_{12}$, on a
\begin{equation}
\bar\psi_2-\bar\psi_1=\psi_2-\psi_1-\alpha-\de((\rho_1+\rho_2)\beta)=
\omega-\alpha-\de\beta=0.
\end{equation}
Les formes $(\bar\psi_i)_{i=1,2}$ coïncident sur l'intersection $\Omega_{12}$,
ce sont donc les restrictions d'une forme globale sur $\Omega$ qu'on notera 
$\psi$, et qui vérifie $\de\psi=\phi$. Pour minorer le spectre de $\Omega$, 
on doit majorer la norme de $\psi$ en fonction de celle de $\phi$, 
indépendamment du choix de $h_\varepsilon$. On commence par écrire
\begin{eqnarray}\label{pvp:ineg1}
\|\psi\|^2&\leq&\|\bar\psi_1\|^2+\|\bar\psi_2\|^2\leq
\|\psi_1+\de(\rho_2\beta)\|^2+\|\psi_2-\alpha-\de(\rho_1\beta)\|^2\nonumber\\
&\leq&3\left(\|\psi_1\|^2+\|\psi_2\|^2+\|\de(\rho_1\beta)\|^2+
\|\de(\rho_2\beta)\|^2+\|\alpha\|^2\right).
\end{eqnarray}
 Dans l'inégalité précédente, chaque norme considérée est relative
au domaine $\Omega_i$ sur lequel la forme est définie. On peut
majorer les premiers termes en utilisant le fait que $\mu_{p,1}(\Omega_i)
\|\psi_i\|^2\leq\|\phi_i\|^2\leq\|\phi\|^2$. Pour les termes faisant
intervenir $\alpha$ et $\beta$, on se ramène aux normes sur $\Omega_{12}$:
comme les supports des formes $\rho_i\beta$ sont contenus dans $\Omega_{12}$
leur norme ne change pas quand on se restreint à $\Omega_{12}$, et
en utilisant le le fait que $\alpha$ est de degré $p$ inférieur à $k$
et qu'elle est de longueur constante pour la métrique $g$ on a 
\begin{equation}
\|\alpha\|_{\Omega_2,h_\varepsilon^2g}\leq\|\alpha\|_{\Omega_2,g}
\leq \frac{\Vol(\Omega_2,g)}{\Vol(\Omega_{12})}\|\alpha\|_{\Omega_{12}}.
\end{equation}
On peut alors majorer les termes restants de l'inégalité (\ref{pvp:ineg1})
en commençant par
\begin{equation}
\frac12\|\de(\rho_i\beta)\|^2\leq c_\rho\|\beta\|^2+\|\de\beta\|^2
\leq\|\de\beta\|^2\left(\frac{c_\rho}{\mu_{p-1,1}(\Omega_{12})}+1
\right),
\end{equation}
où $c_\rho=\sup_i\|\nabla\rho_i\|_\infty^2$.
On majore les normes de $\alpha$ et $\de\beta$ en partant de l'égalité
$\alpha+\de\beta=\psi_2-\psi_1$ et en utilisant le fait que $\alpha$
est harmonique avec condition de bord absolue, donc orthogonale 
aux formes exactes, ce qui donne
\begin{eqnarray}
\|\de\beta\|^2&\leq&\|\psi_2-\psi_1\|^2\leq2\|\psi_1\|^2+2\|\psi_2\|^2
\nonumber\\ &\leq&
\left(\frac2{\mu_{p,1}(\Omega_1)}+\frac2{\mu_{p,1}(\Omega_2)}\right)\|\phi\|^2
\end{eqnarray}
et de la même manière,
\begin{equation}
\|\alpha\|^2\leq\left(\frac2{\mu_{p,1}(\Omega_1)}+\frac2{\mu_{p,1}(\Omega_2)}
\right)\|\phi\|^2.
\end{equation}
On a donc une majoration du quotient $\|\psi\|^2/\|\phi\|^2$ en fonction
de $c_\rho$, de $\mu_{p-1,1}(\Omega_{12})$, des 
$(\mu_{p,1}(\Omega_i))_{i=1,2}$ et du rapport $\frac{\Vol(\Omega_2,g)}
{\Vol(\Omega_{12})}$. Seuls $\mu_{p,1}(\Omega_2)$ et $\frac{\Vol(\Omega_2,g)}
{\Vol(\Omega_{12})}$ dépendent du choix de $h_\varepsilon$, et ils sont 
uniformément minorés.
\end{demo} 

On a finalement une minoration uniforme de $\mu_{p,2}(M,g_\varepsilon)$
et $\mu_{q,1}(M,g_\varepsilon)$ pour $1\leq q\leq k+1$, $q\neq p$, 
l'application du lemme~\ref{pvp:mg} permet donc de conclure la démonstration 
du lemme~\ref{intro:lem}.

\subsection{Propriétés de convergence et de stabilité du spectre}
Pour prescrire le spectre dans classe conforme, nous allons utiliser
des techniques déjà mises en \oe uvre par P.~Guérini dans \cite{gu04} et
que nous allons rappeler ici. Nous expliquerons ensuite comment la
construction géométrique considérée peut être réalisée de manière conforme
à partir d'une métrique quelconque.

\begin{figure}[h]
\begin{center}
\begin{picture}(0,0)%
\includegraphics{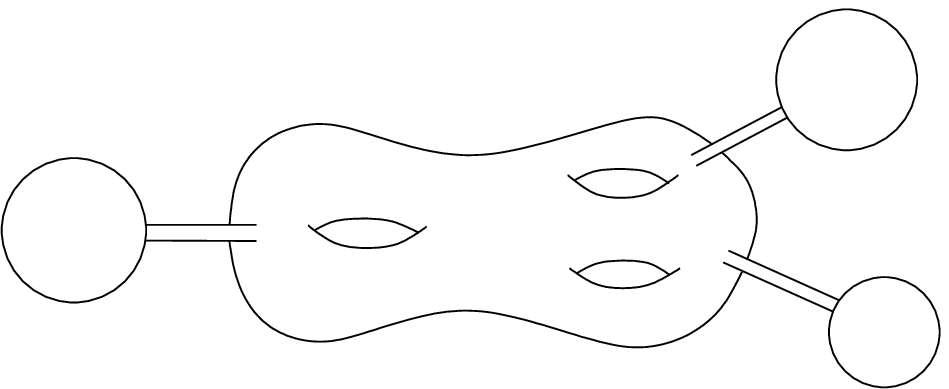}%
\end{picture}%
\setlength{\unitlength}{4144sp}%
\begingroup\makeatletter\ifx\SetFigFontNFSS\undefined%
\gdef\SetFigFontNFSS#1#2#3#4#5{%
  \reset@font\fontsize{#1}{#2pt}%
  \fontfamily{#3}\fontseries{#4}\fontshape{#5}%
  \selectfont}%
\fi\endgroup%
\begin{picture}(4422,1742)(109,-1429)
\put(2131,-136){\makebox(0,0)[lb]{\smash{{\SetFigFontNFSS{12}{14.4}{\rmdefault}{\mddefault}{\updefault}{\color[rgb]{0,0,0}$M_1$}%
}}}}
\put(4411,-16){\makebox(0,0)[lb]{\smash{{\SetFigFontNFSS{12}{14.4}{\rmdefault}{\mddefault}{\updefault}{\color[rgb]{0,0,0}$M_3$}%
}}}}
\put(4516,-1201){\makebox(0,0)[lb]{\smash{{\SetFigFontNFSS{12}{14.4}{\rmdefault}{\mddefault}{\updefault}{\color[rgb]{0,0,0}$M_4$}%
}}}}
\put(316,-211){\makebox(0,0)[lb]{\smash{{\SetFigFontNFSS{12}{14.4}{\rmdefault}{\mddefault}{\updefault}{\color[rgb]{0,0,0}$M_2$}%
}}}}
\end{picture}%
\end{center}
\caption{\label{global:anses}}
\end{figure}

 Un premier outil est le résultat de convergence de spectre obtenu
par C.~Anné et B.~Colbois dans \cite{ac95} pour les variétés compactes
reliées par des anses fines. considérons une famille finie de variétés
compactes $(M_j,g_j)_{i=1}^K$ qu'on relie entre elles par des anses fines,
isométriques au produit d'une sphère $(S^{n-1},\varepsilon^2 
g_{\mathrm{can}})$ 
par un intervalle (voir figure \ref{global:anses}). En notant 
$(\tilde M,g_\varepsilon)$ la variété obtenue, qui est difféomorphe
à $M_1\#M_2\#\ldots\#M_K$, on a alors :
\begin{theo}\label{global:th1}
Si, pour $p\in\{1,\cdots,n-1\}$, on note $\mu'_{p,1}\leq\mu'_{p,2}\leq\ldots$
la réunion des spectres $(\mu_{p,i}(M_j,g_j))_{i,j}$, on a pour tout 
$i\in\N^*$
$$\lim_{\varepsilon\to0}\mu_{p,i}(\tilde M,g_\varepsilon)=\mu'_{p,i}.$$ 
\end{theo}

Décrivons rapidement comment prescrire le spectre (les détails seront
précisés au paragraphe suivant): on commence par fixer un réel strictement 
positif $\delta$ tel que
\begin{equation}
\delta<\inf_{p<k,i<N}\left\{\frac{\nu_{p,i+1}-\nu_{p,i}}2\right\}
\textrm{ et } \delta<V_0
\end{equation}
et choisir un réel $V\in[V_0-\delta,V_0+\delta]$.

Pour tout $1\leq p\leq k$ et tout $1\leq i\leq N$, on se donne un réel 
$\xi_{p,i}\in[\nu_{p,i}-\delta,\nu_{p,i}+\delta]$ et une métrique 
$g_{p,i}$ sur la sphère $S^n$ telle que $\Vol(S^n,g_{p,i})<\frac{V}{Nk}$ et
$\mu_{p,1}(S^n,g_{p,i})=\xi_{p,i}$, toutes les autres valeurs propres de 
la sphère étant supérieures à $\sup_{p,i}\{\nu_{p,i}\}+\delta$.
On munit $M$ d'une métrique telle que $\mu_{p,1}(M)>\sup_{p,i}\{\nu_{p,i}\}+
\delta$ pour tout $p$ compris entre 1 et $n-2$, et $\Vol(M)=
V-\sum_{p,i}\Vol(S^n,g_{p,i})$.  Si on attache les sphères à $M$
par des anses fines comme dans la figure \ref{global:anses}, 
on obtient une variété difféomorphe à $M$,
et le théorème \ref{global:th1} nous
donne alors une famille de métrique $(g_\varepsilon)$ sur $M$ telle 
que $\Vol(M,g_\varepsilon)$ tend vers $V$ et $\mu_{p,i}(M,g_\varepsilon)$
tend vers $\xi_{p,i}$ pour tout $p$ et $i$ quand $\varepsilon$ tend
vers zéro.
Les métriques $g_\varepsilon$ sont singulières sur le bord des anses,
mais on verra qu'on peut les lisser sans perturber la convergence du
spectre.

On peut alors utiliser l'argument de stabilité développé par 
Y.~Colin~de~Verdière dans \cite{cdv86}. On s'appuie sur le
\begin{lem}\label{global:lem}
Soit $(\Phi_\varepsilon)_\varepsilon$ une famille d'applications continues
d'une boule fermée $B_0\subset\R^m$  dans $\R^m$ qui converge uniformément 
vers l'identité quand $\varepsilon$ tend vers $0$, et $x_0$ un point 
intérieur à $B_0$. 

Si $\varepsilon$ est suffisamment petit, alors $x_0$ est contenu dans
l'image de $\Phi_\varepsilon$.
\end{lem}
On applique ce lemme à l'espace $\R^{kN+1}$ avec 
\begin{equation}\label{global:phi}
\Phi_\varepsilon:\left\{\begin{array}{ccl}
[V_0-\delta,V_0+\delta]\times\displaystyle\prod_{p,i} 
[\nu_{p,i}-\delta,\nu_{p,i}+\delta]&\to&\R^{kN+1}\\
(V,\xi_{1,1},\ldots,\xi_{1,N},\ldots,\xi_{k,1},\ldots,\xi_{k,N})&\mapsto&
(\Vol(M,g_\varepsilon),\\
&&\mu_{p,i}(M,g_\varepsilon))\\
\end{array}\right.
\end{equation}
et $x_0=(V_0,\nu_{1,1}\ldots,\nu_{1,N},\nu_{k,1},\ldots,\nu_{k,N})$.
Selon le théorème \ref{global:th1}, $\Phi_\varepsilon$ converge
simplement vers l'identité quand $\varepsilon\to0$. Mais comme
ces applications sont continues et que $B_0$ est compact, le
théorème de Dini garantit que la convergence est uniforme. On peut
donc trouver un jeu de paramètres $(V,\xi_{i,j})$ et un $\varepsilon$
telle que la métrique $g_\varepsilon$ vérifie les conclusions du théorème.
 Comme
$\mu_{k+1,1}(M,g_\varepsilon)$ tend vers une valeur supérieure à
$\sup_{p,i}\{\nu_{p,i}\}+\delta$, on peut choisir $\varepsilon$ 
suffisamment petit pour que la condition $\mu_{k+1,1}(M,g)>
\sup_{p,i}\{\nu_{p,i}\}$ soit vérifiée.

Notons qu'il est est essentiel que pour chaque valeur de $p$ les 
termes de la suite $(\nu_{p,i})_{i=1}^N$ soient distincts, afin qu'on 
puisse choisir $B_0$ et 
$x_0$ tels que $x_0$ soit à l'intérieur de $B_0$.

\subsection{Prescription conforme}

Il reste à montrer que la construction précédente peut être réalisée
de manière conforme. 

La première étape consiste, pour chaque $\xi_{p_0,i}$, à trouver
la métrique $g_{p_0,i}$ correspondante sur la sphère. Pour cela, on
applique le lemme \ref{intro:lem} en considérant la classe conforme 
de la métrique canonique, et
\begin{equation}
\varepsilon<\frac{c\cdot\xi_{p_0,i}}{\sup_{p,i}\{\nu_{p,i}\}+\delta}.
\end{equation}
On obtient une métrique $\bar g_{p_0,i}$ telle que $\mu_{p_0,1}
(S^n,\bar g_{p_0,i})<\varepsilon$, les autres valeurs propres
étant supérieures à $c$. En posant 
\begin{equation}
g_{p_0,i}=\frac{\mu_{p_0,1}(S^n,\bar g_{p_0,i})}{\xi_{p_0,i}}\bar g_{p_0,i},
\end{equation}
on a $\mu_{p_0,1}(S^n,g_{p_0,i})=\xi_{p_0,i}$, et les autres valeurs
propres sont supérieures à $c\cdot\frac{\xi_{p_0,i}}
{\mu_{p_0,1}(S^n,\bar g_{p_0,i})}>c\cdot\frac{\xi_{p_0,i}}\varepsilon>
\sup_{p,i}\{\nu_{p,i}\}+\delta$. Le lemme \ref{intro:lem} permet en outre
de majorer le volume de la sphère pour la métrique $\bar g_{p_0,i}$,
et donc pour la métrique $g_{p_0,i}$. On peut donc choisir
$g_{p_0,i}$ telle que $\Vol(S^n,g_{p_0,i})<\frac{V}{Nk}$.

On veut ensuite munir $M$ d'une
métrique pour laquelle toutes les valeurs propres sont grandes,
le volume et la classe conforme étant fixés. C'est possible en vertu d'un 
résultat de B.~Colbois et A.~El~Soufi:
\begin{theo}[\cite{ces06}] Si $M$ est une variété riemannienne compacte
de dimension $n\geq4$ et $C$ une classe conforme sur $M$, alors
\begin{equation}
\sup_{g\in C}\inf_{0<p<n-2}\mu_{p,1}(M,g)\Vol(M,g)^{\frac2n}=+\infty.
\end{equation}
\end{theo}
On peut donc se donner une métrique $g$ sur $M$ telle que
\begin{equation}
\Vol(M,g)=V-\sum_{p_0,i}\Vol(S^n,g_{p_0,i})
\end{equation}
et \begin{equation}\mu_{p_0,1}(M,g)>\sup_{p,i}\{\nu_{p,i}\}+\delta
\end{equation}
pour tout $0<p_0<n-1$.

\begin{figure}[h]
\begin{center}
\begin{picture}(0,0)%
\includegraphics{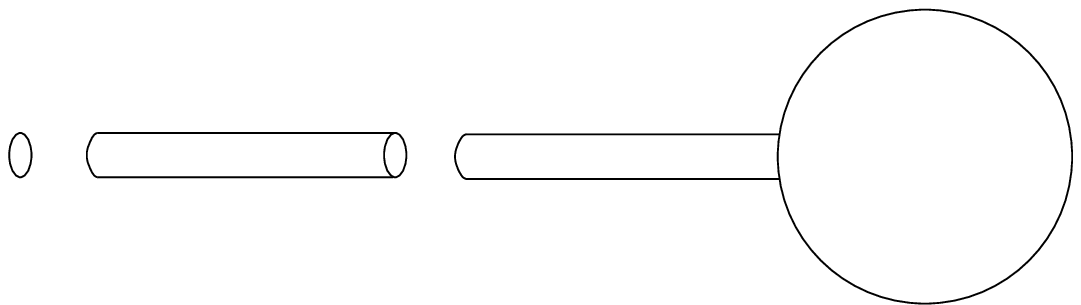}%
\end{picture}%
\setlength{\unitlength}{4144sp}%
\begingroup\makeatletter\ifx\SetFigFontNFSS\undefined%
\gdef\SetFigFontNFSS#1#2#3#4#5{%
  \reset@font\fontsize{#1}{#2pt}%
  \fontfamily{#3}\fontseries{#4}\fontshape{#5}%
  \selectfont}%
\fi\endgroup%
\begin{picture}(5198,1358)(3046,-827)
\put(3061,-556){\makebox(0,0)[lb]{\smash{{\SetFigFontNFSS{10}{12.0}{\rmdefault}{\mddefault}{\updefault}{\color[rgb]{0,0,0}$(B(\varepsilon),g_\varepsilon)$}%
}}}}
\put(3916,-556){\makebox(0,0)[lb]{\smash{{\SetFigFontNFSS{10}{12.0}{\rmdefault}{\mddefault}{\updefault}{\color[rgb]{0,0,0}$(B(\varepsilon),h_{1,\varepsilon}^2g_\varepsilon)$}%
}}}}
\put(5806,-556){\makebox(0,0)[lb]{\smash{{\SetFigFontNFSS{10}{12.0}{\rmdefault}{\mddefault}{\updefault}{\color[rgb]{0,0,0}$(B(\varepsilon),h_{2,\varepsilon}^2g_\varepsilon)$}%
}}}}
\end{picture}%
\end{center}
\caption{\label{global:sphere}}
\end{figure}

On veut maintenant attacher des sphères à $M$ par des anses fines. Un
réel $\varepsilon>0$ petit étant donné, on va ---~temporairement~---
déformer la métrique $g$ de manière non conforme en une métrique
$g_\varepsilon$ telle que cette métrique soit euclidienne sur $kN$ boules
disjointes de rayon $\varepsilon$. On peut choisir $g_\varepsilon$ telle
que $\frac1{\tau(\varepsilon)}g\leq g_\varepsilon\leq \tau(\varepsilon)g$,
avec $\tau(\varepsilon)\to1$ quand $\varepsilon\to0$. Comme l'ont
remarqué B.~Colbois et A.~El~Soufi dans \cite{ces03} et \cite{ces06},
une boule euclidienne $(B(\varepsilon),g_\textrm{euc})$ de rayon $\varepsilon$
peut être déformée de manière conforme en la réunion d'un cylindre 
de rayon $\varepsilon$ et de longueur quelconque, et d'une sphère
homothétique à la sphère canonique privée d'une boule de rayon $\varepsilon$:
si on note $r$ le coordonnée radiale sur $B(\varepsilon)$ et qu'on
définit la fonction $h_{1,\varepsilon}$ par
\begin{equation}
h_{1,\varepsilon}(r)=\left\{\begin{array}{ll}
\frac\varepsilon r&\textrm{si }\varepsilon e^{-\frac L\varepsilon}\leq
r\leq\varepsilon,\\
e^{\frac L\varepsilon}&\textrm{si }0\leq r\leq
\varepsilon e^{-\frac L\varepsilon},\\
\end{array}\right.
\end{equation}
Une fois munie de la métrique $h_{1,\varepsilon}^2g_\textrm{euc}$,
la partie de la boule $B(\varepsilon)$ correspondant à 
$r\in[\varepsilon e^{-\frac L\varepsilon},\varepsilon]$ est isométrique à 
un cylindre de longueur $L$ et de rayon $\varepsilon$ tandis 
que la partie correspondant à $r\in[0,\varepsilon e^{-\frac L\varepsilon}]$ 
est isométrique à une boule euclidienne de rayon $\varepsilon$. Comme cette 
boule peut être projetée stéréographiquement ---~donc de manière conforme~--- 
sur une calotte sphérique quelconque, il existe une fonction 
$h_{2,\varepsilon}$ telle que 
$(B(\varepsilon),h_{2,\varepsilon}^2g_\textrm{euc})$ soit
la réunion d'un cylindre de longueur $L$ et d'un calotte sphérique
dont le bord s'appuie sur le bord du cylindre (voir figure 
\ref{global:sphere}). De plus, on peut choisir
$h_{2,\varepsilon}$ de sorte que lorsque $\varepsilon$ tend vers zéro, la
métrique de la sphère portant la calotte soit fixée. En appliquant cette
déformation sur chacune des $kN$ boules euclidiennes contenues dans 
$(M,g_\varepsilon)$, la variété est isométrique à celle obtenue en
attachant $kN$ sphères à $(M,g_\varepsilon)$ par des anses de rayon 
$\varepsilon$, on se trouve bien dans les condition d'application du
théorème \ref{global:th1}.

Il reste à transplanter sur chacune des sphères la métrique $g_{p_0,i}$
qu'on a défini précédemment. Pour cela, remarquons que dans la construction
de la métrique $g_{p_0,i}$ à partir de la métrique canonique, il y a un
ouvert $\mathcal U$ de la sphère sur lequel la déformation conforme est 
une simple homothétie, dont on notera $\rho$ le rapport (voir
remarque \ref{pvp:rq1}). On fixe un
point $x\in\mathcal U$, et on identifie la calotte sphérique qu'on a
construit à la calotte de $(S^n,\rho^2g_\textrm{can})$ dont le bord
est centré en $x$. On peut alors déformer la métrique sur la calotte
de sorte qu'elle soit isométrique à la restriction de $g_{p_0,i}$.
Si $\varepsilon$ est suffisamment petit, le bord de la calotte est
entièrement contenu dans $\mathcal U$, il n'est donc pas déformé et
la métrique globale est bien continue le long du bord.

En appliquant ces déformations sur chacune des $kN$ boules de 
$(M,g_\varepsilon)$, on construit une fonction $h_{3,\varepsilon}$ sur $M$ 
telle que $(M,h_{3,\varepsilon}^2g_\varepsilon)$ soit formé de
 $(M,g_\varepsilon)$ reliés aux $kN$ sphères $(S^n,g_{p_0,i})$ par des 
anses dont le rayon tend vers zéro quand $\varepsilon$ tend vers zéro. 
Le théorème \ref{global:th1} assure alors que 
\begin{equation}
\mu_{p,i}(M,h_{3,\varepsilon}^2g_\varepsilon)\to\xi_{p,i}\textrm{ quand }
\varepsilon\to0
\end{equation}
pour $1\leq p\leq k$ et $1\leq i\leq N$, les autres valeurs propres
tendant vers des valeurs supérieures à $\sup_{p,i}\{\nu_{p,i}\}+\delta$
et le volume tendant vers $V$. 

La famille de métrique $h_{3,\varepsilon}^2g_\varepsilon$ n'est pas
conforme à la métrique initiale $g$. Cependant, comme on a
$\frac1{\tau(\varepsilon)}g_\varepsilon\leq g\leq 
\tau(\varepsilon)g_\varepsilon$,
on a aussi
\begin{equation}
\frac1{\tau(\varepsilon)}h_{3,\varepsilon}^2g_\varepsilon\leq 
h_{3,\varepsilon}^2g\leq \tau(\varepsilon)h_{3,\varepsilon}^2g_\varepsilon
\end{equation}
avec $\tau(\varepsilon)\to1$.
On peut en outre approcher la famille de fonction 
$(h_{3,\varepsilon})_\varepsilon$ par une famille 
$(h_{4,\varepsilon})_\varepsilon$ formée de fonctions lisses et telle que
$e^{-\varepsilon}<h_{4,\varepsilon}/h_{3,\varepsilon}<e^\varepsilon$. 
La famille $(h_{4,\varepsilon}g)$ est alors constituée de métriques lisses
conformes à $g$ et vérifie 
\begin{equation}
\frac1{e^\varepsilon\tau(\varepsilon)}h_{3,\varepsilon}^2g_\varepsilon\leq
h_{4,\varepsilon}^2g\leq e^\varepsilon\tau(\varepsilon)h_{3,\varepsilon}^2
g_\varepsilon.
\end{equation}
Le lemme~\ref{pvp:lem} assure alors que
\begin{equation}
\lim_{\varepsilon\to0}\mu_{p,i}(M,h_{4,\varepsilon}^2g)=
\lim_{\varepsilon\to0}\mu_{p,i}(M,h_{3,\varepsilon}^2g_\varepsilon)
=\xi_{p,i}
\end{equation}
pour $1\leq p\leq k$ et $1\leq i\leq N$ et que les autres valeurs propres
restent supérieures à $\sup_{p,i}\{\nu_{p,i}\}$ pour $\varepsilon$
suffisamment petit, et en outre $\Vol(M,h_{4,\varepsilon}^2g)$ tend vers
$V$. Cette construction peut être réalisée pour n'importe quelle famille
$(\xi_{p,i})$ du domaine de la fonction $\Phi_\varepsilon$ définie en
(\ref{global:phi}). On est donc en mesure d'appliquer le lemme de 
stabilité~\ref{global:lem} et d'en déduire le théorème~\ref{intro:th}.

\noindent Pierre \textsc{Jammes}\\
Université d'Avignon\\
laboratoire de mathématiques\\
33 rue Louis Pasteur\\
F-84000 Avignon\\
\texttt{Pierre.Jammes@univ-avignon.fr}

\begin{thebibliography}{{Mc}93}
{\small
\makeatletter
\ifx\fonteauteurs\@undefined
\def\fonteauteurs{\scshape}\fi
\makeatother

\bibitem[AC95]{ac95}
\bgroup\fonteauteurs C.~Anné\egroup{} et \bgroup\fonteauteurs
  B.~Colbois\egroup{} -- «~Spectre du laplacien agissant sur les $p$-formes
  différentielles et écrasement d'anses~», {\em Math. Ann.}, 303 (3),
  p.~545--573, 1995.

\bibitem[Am03]{am03}
\bgroup\fonteauteurs B.~Ammann\egroup{} -- «~A spin-conformal lower bound of
  the first positive {D}irac eigenvalue~», {\em Differ. Geom. Appl.}, 18 (1),
  p.~21--32, 2003.

\bibitem[CdV86]{cdv86}
\bgroup\fonteauteurs Y.~Colin~de Verdière\egroup{} -- «~Sur la multiplicité de
  la première valeur propre non nulle du laplacien~», {\em Comment. Math.
  Helv.}, 61 (2), p.~254--270, 1986.

\bibitem[CdV87]{cdv87}
\bgroup\fonteauteurs Y.~Colin~de Verdière\egroup{} -- «~Construction de
  laplaciens dont une partie finie du spectre est donnée~», {\em Ann. scient.
  \'Ec. norm. sup.}, 20 (4), p.~99--615, 1987.

\bibitem[CES03]{ces03}
\bgroup\fonteauteurs B.~Colbois\egroup{} et \bgroup\fonteauteurs
  A.~El~Soufi\egroup{} -- «~Extremal eigenvalues of the {L}aplacian in a
  conformal class of metrics: the ``conformal spectrum''~», {\em Ann. Global
  Anal. Geom.}, 23 (4), p.~337--349, 2003, math.DG/0409316.

\bibitem[CES06]{ces06}
\bgroup\fonteauteurs B.~Colbois\egroup{} et \bgroup\fonteauteurs
  A.~El~Soufi\egroup{} -- «~Eigenvalues of the laplacian acting on $p$-forms
  and metric conformal deformations~», {\em Proc. of Am. Math. Soc.}, 134 (3),
  p.~715--721, 2006, math.DG/0409242.

\bibitem[Co04]{co04}
\bgroup\fonteauteurs B.~Colbois\egroup{} -- «~Spectre conforme et métriques
  extrémales~», {\em Sémin. Théor. Spectr. Géom.}, 22, p.~93--101, 2004.

\bibitem[Da05]{da05}
\bgroup\fonteauteurs M.~Dahl\egroup{} -- «~Prescribing eigenvalues of the
  {D}irac operator~», {\em Manuscripta math.}, 118 (2), p.~191--199, 2005,
  math.DG/0311172.

\bibitem[Do82]{do82}
\bgroup\fonteauteurs J.~Dodziuk\egroup{} -- «~Eigenvalues of the {L}aplacian on
  forms~», {\em Proc. of Am. Math. Soc.}, 85, p.~438--443, 1982.

\bibitem[ESI86]{esi86}
\bgroup\fonteauteurs A.~El~Soufi\egroup{} et \bgroup\fonteauteurs
  S.~Ilias\egroup{} -- «~Immersions minimales, première valeur propre du
  laplacien et volume conforme~», {\em Math. Ann.}, 275 (2), p.~257--267, 1986.

\bibitem[GP95]{gp95}
\bgroup\fonteauteurs G.~Gentile\egroup{} et \bgroup\fonteauteurs
  V.~Pagliara\egroup{} -- «~Riemannian metrics with large first eigenvalue on
  forms of degree $p$~», {\em Proc. of Am. Math. Soc.}, 123 (12),
  p.~3855--3858, 1995.

\bibitem[Gu04]{gu04}
\bgroup\fonteauteurs P.~Guérini\egroup{} -- «~Prescription du spectre du
  laplacien de {H}odge-de~{R}ham~», {\em Ann. scient. \'Ec. norm. sup.}, 37
  (2), p.~270--303, 2004.

\bibitem[Ja04]{ja04}
\bgroup\fonteauteurs P.~Jammes\egroup{} -- «~Petites valeurs propres des fibrés
  principaux en tores~», prépublication, 2004, math.DG/0404536.

\bibitem[Ja06]{ja06}
\bgroup\fonteauteurs P.~Jammes\egroup{} -- «~Minoration conforme du spectre du
  laplacien de {H}odge-de~{R}ham~», prépublication, 2006, math.DG/0604591.

\bibitem[Ko93]{ko93}
\bgroup\fonteauteurs N.~Korevaar\egroup{} -- «~Upper bounds for eigenvalues of
  conformal metrics~», {\em J. differ. geom.}, 37 (1), p.~73--93, 1993.

\bibitem[Lo86]{lo86}
\bgroup\fonteauteurs J.~Lott\egroup{} -- «~Eigenvalue bounds for the {D}irac
  operator~», {\em Pacific J. of Math.}, 125 (1), p.~117--126, 1986.

\bibitem[Lo96]{lo96}
\bgroup\fonteauteurs J.~Lohkamp\egroup{} -- «~Discontinuity of geometric
  expansions~», {\em Comment. Math. Helv.}, 71 (2), p.~213--228, 1996.

\bibitem[{Mc}93]{mc93}
\bgroup\fonteauteurs J.~{Mc}{Gowan}\egroup{} -- «~The $p$-spectrum of the
  {L}aplacian on compact hyperbolic three manifolds~», {\em Math. Ann.}, 297
  (4), p.~725--745, 1993.

}\end{thebibliography}
\end{document}